\documentclass[11pt,a4paper]{article}
\usepackage{amssymb,latexsym,a4wide}
\def\exp{\mathop{\rm exp}\nolimits}

\def\proof{{\bf Proof}.\ }

\newtheorem{formula}{}[section]
\newtheorem{proposition}[formula]{Proposition}
\newtheorem{definition}[formula]{Definition}
\newtheorem{corollary}[formula]{Corollary}
\newtheorem{remark}[formula]{Remark}
\newtheorem{lemma}[formula]{Lemma}
\newtheorem{theorem}[formula]{Theorem}

\def\prop{\begin{proposition}}
\def\propl#1{\begin{proposition}\label{#1}}
\def\eprop{\end{proposition}}
\def\exp{\begin{example}}
\def\expl#1{\begin{example}\label{#1}}
\def\eexp{\end{example}}
\def\thrm{\begin{theorem}}
\def\thrml#1{\begin{theorem}\label{#1}}
\def\ethrm{\end{theorem}}
\def\rmrk{\begin{remark}}
\def\rmrkl#1{\begin{remark}\label{#1}}
\def\ermrk{\end{remark}}
\def\dfntn{\begin{definition}}
\def\dfntnl#1{\begin{definition}\label{#1}}
\def\edfntn{\end{definition}}
\def\nmrt{\begin{enumerate}}
\def\enmrt{\end{enumerate}}

\def\qtn{\begin{equation}}
\def\qtnl#1{\begin{equation}\label{#1}}
\def\eqtn{\end{equation}}
\def\lmm{\begin{lemma}}
\def\lmml#1{\begin{lemma}\label{#1}}
\def\elmm{\end{lemma}}
\def\crllr{\begin{corollary}}
\def\crllrl#1{\begin{corollary}\label{#1}}
\def\ecrllr{\end{corollary}}
\begin{document}
\title{On the complexity of solving ordinary differential equations in terms of Puiseux series}
\author{
Ali AYAD \\[-1pt]
\small Campus scientifique de Beaulieu, IRMAR \\[-3pt]
\small Universit\'e Rennes 1, 35042, Rennes, France\\[-3pt]
{\tt \small ali.ayad@univ-rennes1.fr}\\[-3pt]
}
\date{}

\maketitle

\begin{abstract}
We prove that the binary complexity of solving ordinary polynomial differential equations in terms of Puiseux series is single exponential in the number of terms in the series. Such a bound was given by Grigoriev~\cite{Gri9} for {\it Riccatti} differential polynomials associated to ordinary linear differential operators. In this paper, we get the same bound for arbitrary differential polynomials. The algorithm is based on a differential version of the Newton-Puiseux procedure for algebraic equations.

\end{abstract}

\section*{Introduction}

In this paper we are interested in solving ordinary polynomial differential equations. For such equations we are looking to compute solutions in the set of formal power series and more generality in the set of Puiseux series. There is no algorithm which decide whether a polynomial differential equation has Puiseux series as solutions. We get an algorithm which computes a finite extension of the ground field which generates the coefficients of the solutions. Algorithms which estimate the coefficients of the solutions are given in~\cite{Mah, Was}. 
\newline
\newline
In order to analyse the binary complexity of factoring ordinary linear differential operators, Grigoriev~\cite{Gri9} describes an algorithm which computes a fundamental system of solutions of the {\it Riccatti} equation associated to an ordinary linear differential operator. The binary complexity of this algorithm is single exponential in the order $n$ of the linear differential operator. There are also algorithms for computing series solutions with real exponents~\cite{GriSin, Can1, Del, Can} and complex exponents~\cite{Del}.
\newline 
\newline
Let $K=\mathbb{Q}(T_1, \dots ,T_l)[\eta]$ be a finite extension of a finitely generated field over $\mathbb{Q}$. The variables $T_1, \dots ,T_l$ are algebraically independent over $\mathbb{Q}$ and $\eta$ is an algebraic element over the field $\mathbb{Q}(T_1, \dots ,T_l)$ with the minimal polynomial $\phi \in \mathbb{Z}[T_1, \dots ,T_l][Z]$. Let $\overline{K}$ be an algebraic closure of $K$ and consider the two fields:
$$L= \cup_{\nu \in \mathbb{N}^*}K((x^{\frac{1}{\nu}})), \quad \mathcal{L}= \cup_{\nu \in \mathbb{N}^*}\overline{K}((x^{\frac{1}{\nu}}))$$ 
which are the fields of fraction-power series of $x$ over $K$ (respectively $\overline{K}$), i.e., the fields of Puiseux series of $x$ with coefficients in $K$ (respectively $\overline{K}$). Each element $\psi \in L$ (respectively $\psi \in \mathcal{L}$) can be represented in the form $\psi =\sum_{i\in \mathbb{Q}}c_ix^i, \quad c_i\in K$ (respectively $c_i\in \overline{K}$). The order of $\psi$ is defined by $ord (\psi):= \min \{i\in \mathbb{Q}, c_i\neq 0\}$. The fields $L$ and $\mathcal{L}$ are differential fields with the differential operator 
$$\frac{d}{dx}(\psi)=\sum_{i\in \mathbb{Q}}ic_ix^{i-1}.$$
\newline
Let $F(y_0, \dots , y_n)$ be a polynomial on the variables $y_0, \dots , y_n$ with coefficients in $L$ and consider the associated ordinary differential equation $F(y, \frac{dy}{dx}, \dots , \frac{d^ny}{dx})=0$ which will be denoted by $F(y)=0$. We will describe all the solutions of the differential equation $F(y)=0$ in $\mathcal{L}$ by a differential version of the Newton polygon process. First write $F$ in the form:
$$F=\sum_{i\in \mathbb{Q}, \alpha \in A}f_{i, \alpha}x^iy_0^{\alpha_0} \cdots y_n^{\alpha_n}, \quad f_{i, \alpha}\in K$$ 
where $\alpha=(\alpha_0,\dots ,\alpha_n)$ belongs to a finite subset $A$ of $\mathbb{N}^{n+1}$. The order of $F$ is defined by $ord (F):= \min \, \{i\in \mathbb{Q}; f_{i, \alpha}\neq 0 \, \, \textrm{for a certain}\, \,  \alpha \}$. Without loss of generality we can suppose that each coefficient $f_{i, \alpha}\in \mathbb{Z}[T_1, \dots ,T_l][\eta]$ and so it can be written in the form
$$f_{i, \alpha}=\sum_{j, j_1,\dots ,j_l}b_{j, j_1,\dots ,j_l}T_1^{j_1}\cdots T_l^{j_l}\eta^j, \quad b_{j, j_1,\dots ,j_l}\in \mathbb{Z}.$$ 
We define the degree of $F$ w.r.t. $x$ by $\deg_x(F)=\max \, \{i\in \mathbb{Q}; \, \, f_{i, \alpha}\neq 0 \, \, \textrm{for a certain}\, \,  \alpha\}$ (it can be equal to $+\infty$), the degree of $F$ w.r.t. $T_1, \dots ,T_l$ by $\deg_{T_1, \dots ,T_l}(F)=\max \, \{\deg_{T_1, \dots ,T_l}(f_{i, \alpha}); \, \, i\in \mathbb{Q}, \alpha \in A\}$. We denote by $l(b)$ the binary length of an integer $b$. 
The binary length of $F$ is defined by $l(F)=\max \, \{l(f_{i, \alpha}); \, \, i\in \mathbb{Q}, \alpha \in A\}$ where $l(f_{i, \alpha})$ is the maximum of the binary lengths of its coefficients in $\mathbb{Z}$. 
We can define in the same manner the degrees and the binary length of $\phi$. To estimate the binary complexity of the algorithm of this paper we suppose that $\deg_Z (\phi)\leq d_0, \, \, \deg_{T_1, \dots ,T_l}(\phi)\leq d_1, \, \, l(\phi)\leq M_1, \, \, \deg_{y_0, \dots ,y_n}(F)\leq d,\, \, \deg_{T_1, \dots ,T_l}(F)\leq d_2, \, \, \deg_x(F)\leq d_3$ ($d_3$ can be equal to $+\infty$) and $l(F)\leq M_2$. 
\newline 
\newline 
In the sequel, we will use the following quantities: for each $i\in \mathbb{Q}$, let 
$$R(i)=d_2(dd_0d_1)^{O(il)}$$
and
$$S(i)=n^id_2(dd_0d_1)^{O(il)}(M_1M_2)^{O(i)}\log_2^{i}(dd_3).$$
We can now describe the main theorem of this paper:

\thrml{maintheorem}
Let $F(y)=0$ be a polynomial differential equation with the above bounds. There is an algorithm which computes all Puiseux series solutions of $F(y)=0$ with coefficients in $\overline{K}$, i.e., all solutions of $F(y)=0$ in $\mathcal{L}$. Namely, for each solution $\psi =\sum_{i\in \mathbb{Q}}c_ix^i \in \mathcal{L}$ of $F(y)=0$, the algorithm computes an integer $\nu \in \mathbb{N}^*$ such that for each $i\in \mathbb{Q}$,  it computes a finite extension $K_1=K[\theta]$ of $K$ where $\theta$ is an algebraic element over $K$ computed with its minimal polynomial $\Phi\in K[Z]$ such that $\sum_{j\leq i,\, j \in \mathbb{Q}}c_jx^j \in K_1((x^{\frac{1}{\nu}}))$. Moreover, for any $j\leq i,\, j \in \mathbb{Q}$, we have the following bounds:
\newline 
\newline
- $\deg_Z(\Phi)\leq d^i$.
\newline 
\newline
- $\deg_{T_1, \dots ,T_l}(\Phi), \, \deg_{T_1, \dots ,T_l}(c_j)\leq R(i)$.
\newline 
\newline
- $l(\Phi), \, l(c_j) \leq S(i)$.
\newline 
\newline
- The binary complexity of this computation is $S(i)$.
\ethrm

\rmrkl{}
i) By the corollary of Lemma 3.1 of~\cite{GriSin}, the integer $\nu$ of Theorem~\ref{maintheorem} is constant, i.e., independent of $i$. This constant depends only on the solution $\psi$.
\newline
\newline
ii) In general, we cannot compute a finite extension $K_1$ of $K$ which contains all the coefficients of all the solutions of $F(y)=0$ in $\mathcal{L}$ either an integer $\nu \in \mathbb{N}^*$ such that all the solutions of $F(y)=0$ (in $\mathcal{L}$) are in $K_1((x^{\frac{1}{\nu}}))$. Namely, if we consider the polynomial 
$$F(y_0, y_1, y_2)=xy_0y_2-xy_1^2+y_0y_1$$
then $\psi=cx^{\mu}$ is a solution of $F(y)=0$ in $\mathcal{L}$ for all $c\in \mathbb{C}$ and all $\mu \in \mathbb{Q}$.
\ermrk
In a forthcoming paper we try to get polynomial binary complexity for computing Puiseux series solutions of polynomial differential equations as in the algebraic case~\cite{Chi86, Chi92}. 
\newline
\newline
The paper is organized as follow. In section 1, we give a description of the Newton polygon associated to polynomial differential equations. Operations on these Newton polygons will be discussed in section 2. The algorithm with its binary complexity analysis are described in section 3.

\section{Newton polygons of polynomial differential equations}

Let $F$ be a differential polynomial as in the introduction. We now define the Newton polygon of $F$. For every couple $(i, \alpha)\in \mathbb{Q}\times A$ such that $f_{i, \alpha}\neq 0$ (i.e., every existing term in $F$) we mark the point 
$$P_{i, \alpha}:= (i-\alpha_1-2\alpha_2-\cdots -n\alpha_n, \alpha_0+\alpha_1+ \cdots +\alpha_n)\in \mathbb{Q}\times \mathbb{N}.$$
We denote by $P(F)$ the set of all the points $P_{i, \alpha}$. The convex hull of these points and $(+\infty , 0)$ in the plane $\mathbb{R}^2$ is denoted by $\mathcal{N}(F)$ and is called the Newton polygon of the differential equation $F(y)=0$ in the neighborhood of $x=0$. If $\deg_{y_0,\dots ,y_n}(F)=m$, then $\mathcal{N}(F)$ is situated between the two lines $y=0$ and $y=m$. For each $(a, b)\in \mathbb{Q}^2\setminus \{(0, 0)\}$ we define the set 
$$N(F, a, b):= \{(u, v)\in P(F), \, \, \forall (u', v')\in P(F), \quad au'+bv'\geq au+bv \}.$$
A point $P_{i, \alpha}\in P(F)$ is a vertex of the Newton polygon $\mathcal{N}(F)$ if there exist $(a, b)\in \mathbb{Q}^2\setminus \{(0, 0)\}$ such that $N(F, a, b)=\{P_{i, \alpha}\}$. We remark that $\mathcal{N}(F)$ has a finite number of vertices. A pair of different vertices $e= (P_{i, \alpha}, P_{i', \alpha'})$ forms an edge of $\mathcal{N}(F)$ if there exist $(a, b)\in \mathbb{Q}^2\setminus \{(0, 0)\}$ such that $e\subset N(F, a, b)$. We denote by $E(F)$ (respectively $V(F)$) the set of all the edges $e$ (respectively all the vertices $p$) of $\mathcal{N}(F)$ for which $a>0$ and $b\geq 0$ in the previous definitions. 
It  is easy to prove that if $e\in E(F)$, then there exists a unique pair $(a(e), b(e))\in \mathbb{Z}^2$ such that $GCD(a(e), b(e))=1, \, \, a(e)>0, \, \, b(e)\geq 0$ and $e\subset N(F, a(e), b(e))$. By the inclination of a line we mean the negative inverse of its geometric slope. If $e\in E(F)$, we can prove that the fraction $\mu_e:=\frac{b(e)}{a(e)}\in \mathbb{Q}$ is the inclination of the straight line passing through the edge $e$. If $p\in V(F)$ and $N(F, a, b)=\{p\}$ for a certain $(a, b)$, then the fraction $\mu:=\frac{b}{a}\in \mathbb{Q}$ is the inclination of a straight line which intersects $\mathcal{N}(F)$ exactly in the vertex $p$. 
\newline 
\newline
For each $e\in E(F)$ we define the univariate polynomial (in a new variable $C$)
$$H_{(F, e)}(C):= \sum_{P_{i, \alpha}\in N(F, a(e), b(e))}f_{i, \alpha}C^{\alpha_0+\alpha_1+ \cdots +\alpha_n}(\mu_e)_1^{\alpha_1}\cdots (\mu_e)_n^{\alpha_n}\in K[C],$$
where $(\mu_e)_k:=\mu_e(\mu_e-1)\cdots (\mu_e-k+1)$ for any positive integer $k$. We call $H_{(F, e)}(C)$ the {\it characteristic} polynomial of $F$ associated to the edge $e\in E(F)$. Its degree is at most $m=\deg_{y_0,\dots ,y_n}(F)\leq d$.
\newline 
\newline
If $\psi \in \mathcal{L}$ is a solution of the differential equation $F(y)=0$ such that $ord (\psi)= \mu_e$, i.e., $\psi $ has the form $\psi=\sum_{i\in \mathbb{Q}, \, i\geq \mu_e}c_ix^i$, $c_i\in \overline{K}$, then we have that $H_{(F, e)}(c_{\mu_e})=0$, i.e. $c_{\mu_e}$ is a root of the polynomial $H_{(F, e)}$ in $\overline{K}$. This condition is called a necessary initial condition to have a solution of $F(y)=0$ in the form of $\psi$ (see Lemma 1 of~\cite{Can}). In fact, $H_{(F, e)}(c_{\mu_e})$ equals the coefficient of the lowest term in the expansion of $F(\psi (x))$ with indeterminates $\mu_e$ and $c_{\mu_e}$. Let $A_{(F, e)}:=\{c\in \overline{K}, c\neq 0, H_{(F, e)}(c)=0\}$.
\newline 
\newline
For each $p=(u, v)\in V(F)$, let $\mu_1<\mu_2$ be the inclinations of the adjacent edges at $p$ in $\mathcal{N}(F)$, it is easy to prove that for all rational number $\mu=\frac{b}{a}$, $a\in \mathbb{N}^*,\, \, b\in \mathbb{N}$ such that $N(F, a, b)=\{p\}$, we have $\mu_1<\mu <\mu_2$. We associate to $p$ the polynomial 
$$h_{(F, p)}(\mu):= \sum_{P_{i, \alpha}=p}f_{i, \alpha}(\mu)_1^{\alpha_1}\cdots (\mu)_n^{\alpha_n}\in K[\mu],$$
which is called the {\it indicial} polynomial of $F$ associated to the vertex $p$ (here $\mu$ is considered as an indeterminate). Let $H_{(F, p)}(C)= C^vh_{(F, p)}(\mu)$ defined as above for edges $e\in E(F)$. Let $A_{(F, p)}:=\{\mu \in \mathbb{Q}, \, \,  \mu_1<\mu <\mu_2; \, \, h_{(F, p)}(\mu)=0\}$.

\rmrkl{}
Let $p=(u, v)\in V(F)$ and $e$ be the edge of $\mathcal{N}(F)$ descending from $p$, then $h_{(F, p)}(\mu_e)$ is the coefficient of the monomial $C^v$ in the expansion of the {\it characteristic} polynomial of $F$ associated to $e$.
\ermrk

\section{Some operations on Newton polygons of differential polynomials}

\subsection{Relation between Newton polygons of a differential polynomial and its partial derivatives}

Let $F$ be a differential polynomial as in the introduction. Write $F$ in the form 
$$F=F_0+\cdots +F_d$$
where $F_s=\sum_{i\in \mathbb{Q}, \, |\alpha |=s}f_{i, \alpha}x^iy_0^{\alpha_0}\cdots y_n^{\alpha_n}$ is the homogeneous part of $F$ of degree $s$ with respect to the indeterminates $y_0,\dots ,y_n$, $\alpha=(\alpha_0, \dots ,\alpha_n)\in A$ and $|\alpha |=\alpha_0+\cdots +\alpha_n$ is the norm of $\alpha$. Then the ordinate of any point of $P(F_s)$ is equal to $s$ and  
$$P(F)=\cup_{0\leq s\leq d}P(F_s).$$ 
Let $0\leq j\leq n$. If there exists an integer $k\geq 1$ such that for all $1\leq s\leq d$ (such that $F_s\neq 0$), $D_{s, j}:=\deg_{y_j}(F_s)\geq k$ then we can prove that $P({\partial^k F \over \partial y_j^k})$ is a translation of $P(F)$ defined by the point $(kj, -k)$, i.e., 

$$P\Big ({\partial^k F \over \partial y_j^k}\Big )=P(F)+\{(kj, -k)\}\, \, \textrm{and then}\, \, \mathcal{N}\Big ({\partial^k F \over \partial y_j^k}\Big )=\mathcal{N}(F)+\{(kj, -k)\}.$$
For any $(a, b)\in \mathbb{Q}^2\setminus \{(0, 0)\}$, we have 
$$N\Big ({\partial^k F \over \partial y_j^k}, a, b \Big )= N(F, a, b)+ \{(kj, -k)\}.$$
Thus the edges of $\mathcal{N}\Big ({\partial^k F \over \partial y_j^k}\Big )$ are exactly the translation of those of $\mathcal{N}(F)$. For each $e\in E\Big ({\partial^k F \over \partial y_j^k}\Big )$, its {\it characteristic} polynomial is 
$$H_{\big ({\partial^k F \over \partial y_j^k},\, e\big )}(C)= \sum_{P_{i, \alpha}\in N(F, a(e), b(e))}f_{i, \alpha}C^{\alpha_0+\alpha_1+ \cdots +\alpha_n-k}(\alpha_j)_k(\mu_e)_1^{\alpha_1}\cdots (\mu_e)_j^{\alpha_j-k}\cdots  (\mu_e)_n^{\alpha_n}\in K[C].$$
For each $p\in V\Big ({\partial^k F \over \partial y_j^k}\Big )$, its {\it indicial} polynomial is
$$h_{\big ({\partial^k F \over \partial y_j^k}, \, p\big )}(\mu)=\sum_{P_{i, \alpha}=p}f_{i, \alpha}(\alpha_j)_k(\mu_e)_1^{\alpha_1}\cdots (\mu_e)_j^{\alpha_j-k}\cdots  (\mu_e)_n^{\alpha_n}\in K[\mu].$$
Let $0\leq j_1\neq j_2\leq n$. If there exist integers $k_1, k_2\geq 1$ such that for all $1\leq s\leq d$, $D_{s, j_2}\geq k_2$ and $\deg_{y_{j_1}}\Big ({\partial^{k_2} F_s \over \partial y_{j_2}^{k_2}}\Big )\geq k_1$ then 
$$P\Big ({\partial^{k_1+k_2}F \over \partial y_{j_1}^{k_1}y_{j_2}^{k_2}}\Big )=P(F)+\{(k_1j_1+k_2j_2, -k_1-k_2)\}$$
and then
$$\mathcal{N}\Big ({\partial^{k_1+k_2} F \over \partial y_{j_1}^{k_1}y_{j_2}^{k_2}}\Big )=\mathcal{N}(F)+\{(k_1j_1+k_2j_2, -k_1-k_2)\}.$$
For any $(a, b)\in \mathbb{Q}^2\setminus \{(0, 0)\}$, we have 
$$N\Big ({\partial^{k_1+k_2}F \over \partial y_{j_1}^{k_1}y_{j_2}^{k_2}}\Big )=N(F, a, b)+\{(k_1j_1+k_2j_2, -k_1-k_2)\}.$$
For each $e\in E\Big ({\partial^{k_1+k_2} F \over \partial y_{j_1}^{k_1}y_{j_2}^{k_2}}\Big )$, its {\it characteristic} polynomial is $H_{\big ({\partial^{k_1+k_2} F \over \partial y_{j_1}^{k_1}y_{j_2}^{k_2}},\, e\big )}(C)=$

$$\sum_{P_{i, \alpha}\in N(F, a(e), b(e))}f_{i, \alpha}C^{\alpha_0+\alpha_1+ \cdots +\alpha_n-k_1-k_2}(\alpha_{j_1})_{k_1}(\alpha_{j_2})_{k_2}(\mu_e)_1^{\alpha_1}\cdots (\mu_e)_{j_1}^{\alpha_{j_1}-k_1}\cdots  (\mu_e)_{j_2}^{\alpha_{j_2}-k_2} \cdots (\mu_e)_n^{\alpha_n}.$$

\lmml{}
Let $k\geq 1$ be an integer such that for all $1\leq s\leq d$ (such that $F_s\neq 0$) and for all $0\leq j\leq n$ we have $D_{s, j}\geq k$. For any $e\in E(F)$, the $k$-th derivative of $H_{(F, e)}\in K[C]$ is given by the formula
$$H_{(F, e)}^{(k)}(C)=\sum_{0\leq k_0,\dots , k_n\leq n}(\mu_e)_1^{k_1}\cdots (\mu_e)_n^{k_n}H_{\big ({\partial^k F \over \partial y_0^{k_0}\cdots \partial y_n^{k_n}},\, e\big )}(C).$$
where the sum ranges over all the partitions $(k_0, \dots ,k_n)$ of $k$, i.e., $k_0+\cdots +k_n=k$.
\elmm
\proof By induction on $k$. $\Box$

\subsection{Newton polygons of sums of differential polynomials}

Let $F$ and $G$ be two differential polynomials of degree less or equal than $d$. Write $F=F_0+\cdots +F_d$ and $G=G_0+\cdots +G_d$ where $F_s$ (respectively $G_s$) is the homogeneous part of $F$ (respectively $G$) of degree $s$ with respect to the indeterminates $y_0,\dots ,y_n$. For each $0\leq s\leq d$ such that $F_s+G_s\neq 0$, let $P_s=(s_1, s)$ be the point of the plane defined by $s_1:=\min \{u; \, \, (u, s)\in P(F_s+G_s)\}$. Then it is easy to prove that the convex hull of all the points $P_s$ for all $0\leq s\leq d$ and $(+\infty , 0)$ in the plane $\mathbb{R}^2$ is the Newton polygon of the differential polynomial $F+G$.

\subsection{Newton polygons of evaluations of differential polynomials}

Let $F$ be a differential polynomial as in the introduction, $0\neq c \in \overline{K}, \, \, \mu  \in \mathbb{Q}$ and $G(y)=F(cx^{\mu}+y)$. We will discuss the construction of the Newton polygon of the differential equation $G(y)=0$ for different values of $c$ and $\mu$. 
\newline
\newline
For each differential monomial $m(y)=f_{i, \alpha}x^iy_0^{\alpha_0}\cdots y_n^{\alpha_n}$ of $F$ with corresponding point $p$, compute 
$$m(cx^{\mu}+y)= f_{i, \alpha}x^i(cx^{\mu}+y_0)^{\alpha_0}(c\mu x^{\mu-1}+ y_1)^{\alpha_1}\cdots (c(\mu)_nx^{\mu -n}+ y_n)^{\alpha_n}.$$
Remark that the corresponding points of the differential monomials of $m(cx^{\mu}+y)$ have ordinate less or equal than $s$ and lie in the line passing through $p$ with inclination $\mu$. There are two possibilities for $\mu$:

\lmml{tactac}
If $\mu =\mu_e$ is the inclination of an edge $e\in E(F)$. For any $0\leq s\leq d$, the vertex of $\mathcal{N}(G)$ of ordinate $s$ corresponds to the differential monomial of $G$ with coefficient equals to
$$q_s(c, \mu_e) := \sum_{0\leq k_0\leq \cdots \leq k_n\leq n}\frac{1}{k_0!\cdots k_n!}H_{\big ({\partial^s F \over \partial y_0^{k_0}\cdots \partial y_n^{k_n}},\, e\big )}(c).$$
where the sum ranges over all the partitions $(k_0, \dots ,k_n)$ of $s$. Its x-coordinate is the minimum of the quantities $i+ \mu (\alpha_0+\cdots +\alpha_n-s)-\alpha_1-2\alpha_2-\cdots -n\alpha_n$ for $i\in \mathbb{Q}$ and $\alpha \in A$. 
\newline
If $\mu_1< \mu < \mu_2$ where $\mu_1$ and $\mu_2$ are the inclinations of the two adjacent edges of a vertex $p=(u, v)\in V(F)$, then for any $0\leq s\leq d$, the vertex of $\mathcal{N}(G)$ of ordinate $s$ corresponds to the differential monomial of $G$ with coefficient equals to
$$c^{v-s}\sum_{0\leq k_0\leq \cdots \leq k_n\leq n}\frac{1}{k_0!\cdots k_n!}h_{\big ({\partial^s F \over \partial y_0^{k_0}\cdots \partial y_n^{k_n}},\, p\big )}(\mu).$$
where the sum ranges over all the partitions $(k_0, \dots ,k_n)$ of $s$.
\elmm
\proof
Let $\mu =\mu_e$ for $e\in E(F)$ and compute 
$$G(y)= \sum_{i\in \mathbb{Q}, \alpha \in A}f_{i, \alpha}x^i(cx^{\mu}+y_0)^{\alpha_0}(c\mu x^{\mu-1}+ y_1)^{\alpha_1}\cdots (c(\mu)_nx^{\mu -n}+ y_n)^{\alpha_n}.$$
For each $0\leq s\leq d$, compute $G_s$ the homogeneous part of $G$ of degree $s$ in $y_0, \dots ,y_n$. We remark that for fixed $i$ and $\alpha$, all the differential monomials of $G_s$ have the same corresponding point which is 
$$(i+ \mu (\alpha_0+\cdots +\alpha_n-s)-\alpha_1-2\alpha_2-\cdots -n\alpha_n, \, s).$$
The x-coordinate of the vertex of $\mathcal{N}(G)$ of ordinate $s$ is the minimum of the x-coordinates $i+ \mu (\alpha_0+\cdots +\alpha_n-s)-\alpha_1-2\alpha_2-\cdots -n\alpha_n$ for $i\in \mathbb{Q}$ and $\alpha \in A$. This minimum is realized by the points $P_{i, \alpha}\in N(F, a(e), b(e))$. This proves the lemma taking into account the formula for the {\it characteristic} polynomial of the derivatives of $F$ in the previous subsections. $\Box$
\newline
\newline
The following lemma is a generalization of Lemma 2.2 of ~\cite{Gri9} which deals with the Newton polygon of the {\it Riccatti} equation associated to a linear ordinary differential equation.

\crllrl{bayen}
Let $\mu =\mu_e$ be the inclination of an edge $e\in E(F)$. The edges of $\mathcal{N}(G)$, situated above the edge $e$ are the same as in $\mathcal{N}(F)$. Let an integer $s_1\geq 0$ be such that $q_s(c, \mu_e)=0$ for all $0\leq s < s_1$ and $q_{s_1}(c, \mu_e)\neq 0$ then $\mathcal{N}(G)$ has a vertex of ordinate $s_1$ and it has at least $s_1$ edges with inclination greater than $\mu_e$.
\ecrllr
\proof
Let $p_{s_1}$ be the vertex of $\mathcal{N}(G)$ of ordinate $s_1$ and $x$-coordinate $x_{p_{s_1}}$ the minimum of the values $i+ \mu (\alpha_0+\cdots +\alpha_n-s_1)-\alpha_1-2\alpha_2-\cdots -n\alpha_n$ for $i\in \mathbb{Q}$ and $\alpha \in A$ (by Lemma~\ref{tactac}). We have $q_{s_1-1}(c, \mu_e)=0$, then the $x$-coordinate of the vertex $p_{s_1-1}$ of ordinate $s_1-1$ is strictly less than the minimum of the values $i+ \mu (\alpha_0+\cdots +\alpha_n-s_1+1)-\alpha_1-2\alpha_2-\cdots -n\alpha_n$ for $i\in \mathbb{Q}$ and $\alpha \in A$. Thus the inclination of the edge joining $p_{s_1}$ and $p_{s_1-1}$ is greater than $\mu$. $\Box$

\crllrl{train}
Let $\mu =\mu_e$ be the inclination of an edge $e\in E(F)$. If $H_{(F, e)}(c)=0$ then the intersection point of the straight line passing through $e$ with the x-axis is not a vertex of $\mathcal{N}(G)$ and $\mathcal{N}(G)$ has an edge with inclination greater than $\mu_e$.
\ecrllr
\proof We have $q_0(c, \mu_e)=H_{(F, e)}(c)=0$, then $s_1\geq 1$. This proves the corollary by applying Corollary~\ref{bayen}. $\Box$

\section{Differential version of the Newton-Puiseux algorithm}

We describe now a differential version of the Newton-Puiseux algorithm to give formal Puiseux series solutions of the differential equation $F(y)=0$. The input of the algorithm is a differential polynomial equation $F(y)=0$ with the bounds described in the introduction. The algorithm will construct a tree $\mathcal{T}$ which depends only on $F$ and on the field $K$. The root of $\mathcal{T}$ is denoted by $\tau_0$. For each node $\tau$ of $\mathcal{T}$, it constructs the following elements:
\newline

- The field $K_{\tau}$ which is a finite extension of $K$.

- The primitive element $\theta_{\tau}$ of the extension $K_{\tau}$ of $K$ with its minimal polynomial $\phi_{\tau}\in K[Z]$.

- An element $c_{\tau}\in K_{\tau}$, a number $\mu_{\tau}\in \mathbb{Q}\cup \{-\infty, +\infty\}$ and an element $y_{\tau}=c_{\tau}x^{\mu_{\tau}}+y_{\tau_1}\in K_{\tau}((x^{\frac{1}{\nu (\tau)}}))$ where $\tau$ is a descendant of $\tau_1$ (here $\mu_{\tau}>\mu_{\tau_1}$) and $\nu (\tau)\in \mathbb{N}^*$. 

- The differential polynomial $F_{\tau}(y)=F(y+y_{\tau})$ with coefficients in $K_{\tau}((x^{\frac{1}{\nu (\tau)}}))$.
\newline
\newline
We define the degree of $\tau$ by $\deg (\tau)=\mu_{\tau}\in \mathbb{Q}$, we have $\deg (\tau)=\deg_x(y_{\tau})$ if $c_{\tau}\neq 0$. The level of the node $\tau$, denoted by $lev (\tau)$, is the distance from $\tau_0$ to $\tau$.
\newline
\newline
For the root $\tau_0$ we have $K_{\tau_0}=K,\, \, \theta_{\tau_0}=1, \, \, \phi_{\tau_0}=Z-1, \, \, c_{\tau_0}=y_{\tau_0}=0, \, \, \deg (\tau_0)=\mu_{\tau_0}=-\infty, \nu (\tau_0)=1$ and $F_{\tau_0}(y)=F(y)$.
\newline
\newline
A node $\tau$ of the tree $\mathcal{T}$ is a leaf of $\mathcal{T}$ if for each $e\in E(F_{\tau})$ and for each $p\in V(F_{\tau})$ we have $\mu_e\leq \deg (\tau)$ and $\mu_2\leq \deg (\tau)$ and $y=0$ is a solution of $F_{\tau}(y)=0$, where $\mu_1<\mu_2$ are the inclinations of the adjacent edges at $p$ in $\mathcal{N}(F)$.
\newline
\newline
The algorithm constructs the tree $\mathcal{T}$ by induction on the level of its nodes. We suppose by induction on $i$ that all the nodes of $\mathcal{T}$ of level $\leq i$ are constructed. Denote by $\mathcal{T}_i$  the set of these nodes. At the $(i+1)$-th step of the induction, for each node $\tau$ of level $i$ which is not a leaf of $\mathcal{T}$ we consider the following two sets:
\newline

- $E'(F_{\tau})= \{e\in E(F_{\tau}), \, \, \mu_e>\deg (\tau)\}$ and 

- $V'(F_{\tau})= \{p\in V(F_{\tau}), \, \, \mu_2> \deg (\tau)\}$.
\newline
\newline
For each $e\in E'(F_{\tau})$, compute a factorization of the polynomial $H_{(F_{\tau}, e)}(C)\in K_{\tau}[C]$ into irreducible factors over the field $K_{\tau}=K[\theta_{\tau}]$ in the form
$$H_{(F_{\tau}, e)}(C)=\lambda_e \prod_jH_j^{k_j}$$
where $0\neq \lambda_e\in K_{\tau}, \, \, k_j\in \mathbb{N}^*$ and $H_j\in K_{\tau}[C]$ are monic and irreducibles over $K_{\tau}$. We can do this factorization by the algorithm of~\cite{Chi3, Chi1, Chi, Gri}.  
The elements of the set $A_{(F_{\tau}, e)}$ correspond to the roots of the factors $H_j\neq C$. We consider a root $c_j\in \overline{K}$ for each factor $H_j\neq C$ and we compute a primitive element $\theta_{j, e, \tau}$ of the finite extension $K_{\tau}[c_j]=K[\theta_{\tau}, c_j]$ of $K$ with its minimal polynomial $\phi_{j, e, \tau}\in K[Z]$ using the algorithm of~\cite{Chi3, Chi, Gri}.
\newline
For each root $c_j$ of $H_j\neq C$ we correspond a son $\sigma$ of $\tau$ such that $\theta_{\sigma}=\theta_{j, e, \tau}$, the field $K_{\sigma}=K[\theta_{j, e, \tau}]$ and the minimal polynomial of $\theta_{\sigma}$ over $K$ is $\phi_{\sigma}=\phi_{j, e, \tau}$. Moreover, $c_{\sigma}=c_j, \, \, \mu_{\sigma}=\mu_e, \, \, y_{\sigma}=c_{\sigma}x^{\mu_{\sigma}}+y_{\tau}$ and $F_{\sigma}(y)=F(y+y_{\sigma})$. For $\nu (\sigma)$, we take $\nu (\sigma)=LCM (\nu (\tau), a(e))$ for example.
\newline
\newline
For each $p\in V'(F_{\tau})$, we consider the {\it indicial} polynomial $h_{(F_{\tau}, p)}(\mu) \in K_{\tau}[\mu]$ of $F_{\tau}$ associated to $p$. To each $\mu \in A_{(F_{\tau}, p)}$ such that $\mu > \deg (\tau)$ and $0\neq c\in \overline{K}$ (where $c$ is given by its minimal polynomial over $K$), we correspond a son $\sigma$ of $\tau$ such that $\theta_{\sigma}=c, \, \, c_{\sigma}=c, \, \, \mu_{\sigma}=\mu$. This completes the description of all the sons of the node $\tau$ of the tree $\mathcal{T}$. 

\rmrkl{}
i) If ($E'(F_{\tau})\neq \emptyset$ or $V'(F_{\tau})\neq \emptyset$) and $y=0$ is a solution of $F_{\tau}(y)=0$ then one of the sons of $\tau$ is a leaf $\sigma$ for which $F_{\sigma}=F_{\tau}, \, \, \mu_{\sigma}=+ \infty$ and $c_{\sigma}=0$.
\newline
ii) For any node $\tau$ of $\mathcal{T}$ such that $\deg (\tau)\neq \infty$, if $y=0$ is not a solution of $F_{\tau}(y)=0$ then $E'(F_{\tau})\neq \emptyset$ by Corollary~\ref{train}.
\ermrk

\subsection{Determination of the solutions of $F(y)=0$ in $\mathcal{L}$ by the leaves of $\mathcal{T}$}

Let $\mathcal{U}$ be the set of all the vertices $\tau$ of $\mathcal{T}$ such that either $\deg (\tau)=+\infty$ and for the ancestor $\tau_1$ of $\tau$ it holds $\deg (\tau_1)<+\infty$ or $\deg (\tau)<+\infty$ and $\tau$ is a leaf of $\mathcal{T}$. For each $\tau \in \mathcal{U}$, there exists a sequence $(\tau_i(\tau))_{i\geq 0}$ of vertices of $\mathcal{T}$ such that $\tau_0(\tau)=\tau_0$ and $\tau_{i+1}(\tau)$ is a son of $\tau_i(\tau)$ for all $i\geq 0$. For each $\tau \in \mathcal{U}$, the element 
$$y_{\tau}=\sum_{i\geq 0}c_{\tau_i(\tau)}x^{\mu_{\tau_i(\tau)}}\in K_{\tau}((x^{\frac{1}{\nu (\tau)}}))$$
is a solution of $F(y)=0$. In fact, there are two possibilities to $\tau$: if $\deg (\tau)=+\infty$ then $y=0$ is a solution of $F_{\tau_1}(y)=0$ where $\tau$ is a son of $\tau_1$ and so $y_{\tau_1}$ is a solution of $F(y)=0$. If $\deg (\tau)<+\infty$ and $\tau$ is a leaf of $\mathcal{T}$ then $y=0$ is a solution of $F_{\tau}(y)=F(y+y_{\tau})=0$ and so $y_{\tau}$ is a solution of $F(y)=0$. This defines a bijection between $\mathcal{U}$ and the set of the solutions of $F(y)=0$ in $\mathcal{L}$.

\subsection{Binary complexity analysis of the Newton-Puiseux algorithm}

We begin by estimating the binary complexity of computing all the sons of the root $\tau_0$ of $\mathcal{T}$. For each $e\in E(F)$, we consider the polynomial $H_{(F, e)}(C)\in K[C]$, its degree w.r.t. $C$ (respectively $T_1, \dots ,T_l$) is bounded by $d$ (respectively $d_2$). We have $\mu_e\leq \frac{d_3}{d}$ and its binary length is $l(\mu_e)\leq O(\log_2(dd_3))$ (using the fact that $\mu_e$ is the inclination of the straight line passing through $e$). Then the binary length of $H_{(F, e)}(C)$ is bounded by $M_2+ndO(\log_2(dd_3))$. By the algorithm of~\cite{Chi3, Chi1, Chi, Gri} the binary complexity of factoring $H_{(F, e)}(C)$ into irreducible polynomials over $K$ is 
$$(dd_1d_2)^{O(l)}(nd_0M_1M_2\log_2(dd_3))^{O(1)}.$$
Moreover, each factor $H_j\in K[C]$ of $H_{(F, e)}(C)$ satisfies the following bounds (see Lemma 1.3 of~\cite{Chi}): $\deg_C(H_j)\leq d$, $\deg_{T_1, \dots ,T_l}(H_j)\leq d_2(dd_0d_1)^{O(1)}$ and 
$$l(H_j)\leq nld_2(dd_0d_1)^{O(1)}M_1M_2\log_2(dd_3).$$
By the induction we suppose that the following bounds hold at the $i$-th step of the algorithm for each node $\tau$ of $\mathcal{T}$ of level $i$: 
\newline
\newline
- $\deg_Z(\phi_{\tau})\leq d^i$.
\newline
\newline
- $\deg_{T_1, \dots ,T_l}(\phi_{\tau}), \, \deg_{T_1, \dots ,T_l}(c_{\tau})\leq R(i)$.
\newline
\newline
- $l(\phi_{\tau}), \, l(c_{\tau})\leq S(i)$ where $R(i)$ and $S(i)$ are as in the introduction.
\newline
\newline
- $\mu_{\tau}\leq i(\frac{d_3}{d})$ and then $l(\mu_{\tau})\leq O(\log_2(idd_3))$.
\newline
\newline
Then we have the following bounds for the differential polynomial $F_{\tau}(y)=F(y+y_{\tau})\in K_{\tau}((x^{\frac{1}{\nu (\tau)}}))[y_0, \dots ,y_n]$:
\newline
\newline
- $\deg_{y_0, \dots ,y_n}(F_{\tau})\leq d$.
\newline
\newline
- $\deg_{T_1, \dots ,T_l}(F_{\tau})\leq d_2+d \deg_{T_1, \dots ,T_l}(c_{\tau})\leq R(i)$.
\newline
\newline
- $\deg_x (F_{\tau})\leq d_3+d\mu_{\tau}\leq (i+1)d_3$.
\newline
\newline
- $l(F_{\tau})\leq M_2+d l(c_{\tau})\leq S(i)$.
\newline
\newline
We compute a primitive element $\eta_1$ of the finite extension $K_{\tau}$ over the field $\mathbb{Q}(T_1, \dots ,T_l)$, i.e., $K_{\tau}=K[\theta_{\tau}]=\mathbb{Q}(T_1, \dots ,T_l)[\eta][\theta_{\tau}]= \mathbb{Q}(T_1, \dots ,T_l)[\eta_1]$ by the corollary of Proposition 1.4 of~\cite{Gri9} (see also section 3 of chapter 1 of~\cite{Chi}). Moreover, $\eta_1= \eta+ \gamma \theta_{\tau}$ where $0\leq \gamma \leq [K_{\tau}: \mathbb{Q}(T_1, \dots ,T_l)] = \deg_Z(\phi)\deg_Z(\phi_{\tau})\leq d^id_0$ and we can compute the monic minimal polynomial $\phi_1\in \mathbb{Q}(T_1, \dots ,T_l)[Z]$ of $\eta_1$ which satisfies the following bounds:
\newline
\newline
- $\deg_Z(\phi_1)\leq d^id_0$
\newline
\newline
- $\deg_{T_1, \dots ,T_l}(\phi_1)
\leq (d^id_0)^{O(1)}$
\newline
\newline
- $l(\phi_1)
\leq S(i)$.
\newline
\newline
- This computation can be done with binary complexity
$
S(i).$
\newline
\newline
For each $e\in E'(F_{\tau})$, we consider the polynomial $H_{(F_{\tau}, e)}(C)\in K_{\tau}[C]$, its degree w.r.t. $C$ (respectively $T_1, \dots ,T_l$) is bounded by $d$ (respectively $R(i)$). We have $\mu_e\leq (i+1)(\frac{d_3}{d})$ and its binary length is $l(\mu_e)\leq O\Big (\log_2((i+1)dd_3)\Big )$. Then the binary length of $H_{(F_{\tau}, e)}(C)$ is bounded by 
$$l(F_{\tau})+ndl(\mu_e)\leq S(i).$$
By the algorithm of~\cite{Chi3, Chi1, Chi, Gri} the binary complexity of factoring $H_{(F_{\tau}, e)}(C)$ into irreducible polynomials over $K_{\tau}=\mathbb{Q}(T_1, \dots ,T_l)[\eta_1]$ is 
$S(i).$
Moreover, each factor $H_j\in K_{\tau}[C]$ of $H_{(F_{\tau}, e)}(C)$ satisfies the following bounds:
\newline
\newline
- $\deg_C(H_j)\leq d$
\newline
\newline
- $\deg_{T_1, \dots ,T_l}(H_j)\leq R(i)$.
\newline
\newline
- $l(H_j)\leq S(i).$
\newline
\newline
Let $c_j\in \overline{K}$ be a root of $H_j$. We can compute by the corollary of Proposition 1.4 of~\cite{Gri9} a primitive element $\theta_{j, e, \tau}$ of the finite extension $K_{\tau}[c_j]=K[\theta_{\tau}, c_j]$ of $K$ with its minimal polynomial $\phi_{j, e, \tau}\in K[Z]$. We can express $\theta_{\sigma} =\theta_{j, e, \tau}$ in the form $\theta_{\sigma}=\theta_{\tau}+ \gamma_jc_j$ where $0\leq \gamma_j\leq \deg_Z(\phi_{\tau})\deg_C(H_j)\leq d^{i+1}$ and $c_{\sigma}=c_j$ in the form 
$$c_{\sigma}=\sum_{0\leq t<d^{i+1}}b_t\theta_{\sigma}^t$$
 where $b_t\in K$. Moreover, the following bounds hold:
\newline
\newline
- $\deg_Z(\phi_{\sigma})\leq d^{i+1}$.
\newline
\newline
- $\deg_{T_1, \dots ,T_l}(\phi_{\sigma}), \, \deg_{T_1, \dots ,T_l}(b_t)\leq R(i)$.
\newline
\newline
- $l(\phi_{\sigma}), \, l(b_t)\leq S(i)$.
\newline
\newline
- $\mu_{\sigma}=\mu_e\leq (i+1)(\frac{d_3}{d})$ and then $l(\mu_{\sigma})\leq O\Big (\log_2((i+1)dd_3)\Big )$.
\newline
\newline
- This computation can be done with binary complexity $S(i)$ and thus the total binary complexity of computing all the sons $\sigma$ of $\tau$ is $S(i)$.

\end{document}